\documentclass[english,oneside]{amsart}
\usepackage[T1]{fontenc}
\usepackage[latin1]{inputenc}
\usepackage{amssymb}

\makeatletter
 \theoremstyle{plain}
\newtheorem{thm}{Theorem}[section]
  \theoremstyle{plain}
  \newtheorem{lem}[thm]{Lemma}

\usepackage{babel}
\makeatother
\begin{document}

\title{On Asymptotics of $q$-Gamma Functions  }

\author{Ruiming Zhang}


\email{ruimingzhang@yahoo.com}


\begin{abstract}
In this paper we derive some asymptotic formulas for the $q$-Gamma
function $\Gamma_{q}(z)$ for $q$ tending to $1$.
\end{abstract}

\date{November 30, 2006}

\curraddr{School of Mathematics\\
Guangxi Normal University\\
Guilin City, Guangxi 541004\\
P. R. China.}

\thanks{This work is partially supported by Guangxi Education Grant No. 200607MS136.}

\maketitle

\section{Introduction\label{sec:intro}}

Given complex numbers \begin{equation}
0<q<1,\quad a\in\mathbb{C},\label{eq:1.1}\end{equation}
we define \cite{Andrews4,Gasper,Ismail2,Koekoek} \begin{equation}
(a;q)_{\infty}:=\prod_{k=0}^{\infty}(1-aq^{k}),\label{eq:1.2}\end{equation}
and the $q$-Gamma function\begin{equation}
\Gamma_{q}(z):=\frac{(q;q)_{\infty}}{(q^{z};q)_{\infty}}(1-q)^{1-z}\quad z\in\mathbb{C}.\label{eq:1.3}\end{equation}
The Euler Gamma function $\Gamma(z)$ is defined as \cite{Andrews4,Gasper,Ismail2,Koekoek}
\begin{equation}
\frac{1}{\Gamma(z)}=z\prod_{k=1}^{\infty}\left(1+\frac{z}{k}\right)\left(1+\frac{1}{k}\right)^{-z},\quad z\in\mathbb{C}.\label{eq:1.4}\end{equation}
The Gamma function satisfies the reflection formula 

\begin{equation}
\Gamma(z)\Gamma(1-z)=\frac{\pi}{\sin\pi z},\quad z\in\mathbb{C},\label{eq:1.5}\end{equation}
and the integral representation\begin{equation}
\Gamma(z)=\int_{0}^{\infty}e^{-t}t^{z-1}dt,\quad\Re(z)>0.\label{eq:1.6}\end{equation}
The Gamma function is a very important function in the theory of special
functions, since all the hypergeometric series are defined in terms
of the shifted factorials $(a)_{n}$, which are quotients of two Gamma
functions \begin{equation}
(a)_{n}:=\frac{\Gamma(a+n)}{\Gamma(a)},\quad a\in\mathbb{C},\quad n\in\mathbb{Z}.\label{eq:1.7}\end{equation}
Similarly, the $q$-Gamma function is also very important in the theory
of the basic hypergeometric series, because all the basic hypergeometric
series are defined in terms of the $q$-shifted factorials $(a;q)_{n}$,
which are scaled quotients of the $q$-Gamma funtions\begin{equation}
(q^{\alpha};q)_{n}=\frac{(1-q)^{n}\Gamma_{q}(\alpha+n)}{\Gamma_{q}(\alpha)},\quad\alpha\in\mathbb{C},\quad n\in\mathbb{Z}.\label{eq:1.8}\end{equation}
W. Gosper heuristically argued that \cite{Andrews4,Gasper,Ismail2,Koekoek}\begin{equation}
\lim_{q\to1}\frac{1}{\Gamma_{q}(z)}=\frac{1}{\Gamma(z)},\quad z\in\mathbb{C}.\label{eq:1.9}\end{equation}
For a rigorous proof of the case $z\in\mathbb{R}$, see \cite{Andrews4}.
In this short note we are going to derive some asymptotic formulas
for $\Gamma_{q}(z)$ as $q\to1$ in two different modes. In the first
mode we let $z\to\infty$ and $q\to1$ simultaneously, while in the
second mode we let $q\to1$ for a fixed $z$. 

\begin{lem}
\label{lem:1}Given any complex number $a$, assume that\begin{equation}
0<\frac{\left|a\right|q^{n}}{1-q}<\frac{1}{2}\label{eq:1.10}\end{equation}
 for some positive integer $n$. Then, for any positive integer $K$,
we have \begin{equation}
\frac{(a;q)_{n}}{(a;q)_{\infty}}=\frac{1}{(aq^{n};q)_{\infty}}:=\sum_{k=0}^{K-1}\frac{(aq^{n})^{k}}{(q;q)_{k}}+r_{1}(a,n,K)\label{eq:1.11}\end{equation}
 with\begin{equation}
\left|r_{1}(a,n,K)\right|\le\frac{2\left(\left|a\right|q^{n}\right)^{K}}{(q;q)_{K}},\label{eq:1.12}\end{equation}
 and \begin{equation}
\frac{(a;q)_{\infty}}{(a;q)_{n}}=(aq^{n};q)_{\infty}:=\sum_{k=0}^{K-1}\frac{q^{k(k-1)/2}}{(q;q)_{k}}(-aq^{n})^{k}+r_{2}(a,n,K)\label{eq:1.13}\end{equation}
 with \begin{equation}
\left|r_{2}(a,n,K)\right|\le\frac{2q^{K(K-1)/2}(\left|a\right|q^{n})^{K}}{(q;q)_{K}}.\label{eq:1.14}\end{equation}

\end{lem}
\begin{proof}
From the $q$-binomial theorem \cite{Andrews4,Gasper,Ismail2,Koekoek}\begin{equation}
\frac{(az;q)_{\infty}}{(z;q)_{\infty}}=\sum_{k=0}^{\infty}\frac{(a;q)_{k}}{(q;q)_{k}}z^{k}\quad a,z\in\mathbb{C},\label{eq:1.15}\end{equation}
we obtain \begin{eqnarray*}
r_{1}(a,n,K) & = & \sum_{k=K}^{\infty}\frac{\left(aq^{n}\right)^{k}}{(q;q)_{k}}=\frac{(aq^{n})^{K}}{(q;q)_{K}}\sum_{k=0}^{\infty}\frac{\left(aq^{n}\right)^{k}}{(q^{K+1};q)_{k}}.\end{eqnarray*}
 Since \[
(q^{K+1};q)_{k}\ge(1-q)^{k}\]
 for $k=0,1,...$, thus, \begin{align*}
 & \left|r_{1}(a,n,K)\right|\le\frac{(\left|a\right|q^{n})^{K}}{(q;q)_{K}}\sum_{k=0}^{\infty}\left(\frac{\left|a\right|q^{n}}{1-q}\right)^{k}\le\frac{2(\left|a\right|q^{n})^{K}}{(q;q)_{K}}.\end{align*}
 Apply a limiting case of \eqref{eq:1.15}, \begin{equation}
(z;q)_{\infty}=\sum_{k=0}^{\infty}\frac{q^{k(k-1)/2}}{(q;q)_{k}}(-z)^{k}\quad z\in\mathbb{C},\label{eq:1.16}\end{equation}
 we get\begin{align*}
r_{2}(a,n,K) & =\frac{q^{K(K-1)/2}(-aq^{n})^{K}}{(q;q)_{K}}\sum_{k=0}^{\infty}\frac{(-aq^{n})^{k}q^{k(k+2K-1)/2}}{(q^{K+1};q)_{k}}.\end{align*}
From the inequalities,\[
\frac{1-q^{k}}{1-q}\ge kq^{k-1},\quad\frac{(q^{K+1};q)_{k}}{(1-q)^{k}}\ge k!q^{k(k+2K-1)/2},\quad\mbox{for }k=0,1,\dots\]
we obtain \begin{align*}
 & \left|r_{2}(a,n,K)\right|\le\frac{q^{K(K-1)/2}(\left|a\right|q^{n})^{K}}{(q;q)_{K}}\sum_{k=0}^{\infty}\frac{1}{k!}\left(\frac{\left|a\right|q^{n}}{1-q}\right)^{k}\\
 & \le\frac{q^{K(K-1)/2}(\left|a\right|q^{n})^{K}}{(q;q)_{K}}\exp(1/2)<\frac{2q^{K(K-1)/2}(\left|a\right|q^{n})^{K}}{(q;q)_{K}}.\end{align*}

\end{proof}
The Jacobi theta functions are defined as \begin{align}
\theta_{1}(z;q):=\theta_{1}(v|\tau) & :=-i\sum_{k=-\infty}^{\infty}(-1)^{k}q^{(k+1/2)^{2}}e^{(2k+1)\pi iv},\label{eq:1.17}\\
\theta_{2}(z;q):=\theta_{2}(v|\tau) & :=\sum_{k=-\infty}^{\infty}q^{(k+1/2)^{2}}e^{(2k+1)\pi iv},\label{eq:1.18}\\
\theta_{3}(z;q):=\theta_{3}(v|\tau) & :=\sum_{k=-\infty}^{\infty}q^{k^{2}}e^{2k\pi iv},\label{eq:1.19}\\
\theta_{4}(z;q):=\theta_{4}(v|\tau) & :=\sum_{k=-\infty}^{\infty}(-1)^{k}q^{k^{2}}e^{2k\pi iv},\label{eq:1.20}\end{align}
 where\begin{equation}
z=e^{2\pi iv},\quad q=e^{\pi i\tau},\quad\Im(\tau)>0.\label{eq:1.21}\end{equation}
 The Jacobi's triple product identities are\begin{align}
\theta_{1}(v|\tau) & =2q^{1/4}\sin\pi v(q^{2};q^{2})_{\infty}(q^{2}e^{2\pi iv};q^{2})_{\infty}(q^{2}e^{-2\pi iv};q^{2})_{\infty},\label{eq:1.22}\\
\theta_{2}(v|\tau) & =2q^{1/4}\cos\pi v(q^{2};q^{2})_{\infty}(-q^{2}e^{2\pi iv};q^{2})_{\infty}(-q^{2}e^{-2\pi iv};q^{2})_{\infty},\label{eq:1.23}\\
\theta_{3}(v|\tau) & =(q^{2};q^{2})_{\infty}(-qe^{2\pi iv};q^{2})_{\infty}(-qe^{-2\pi iv};q^{2})_{\infty},\label{eq:1.24}\\
\theta_{4}(v|\tau) & =(q^{2};q^{2})_{\infty}(qe^{2\pi iv};q^{2})_{\infty}(qe^{-2\pi iv};q^{2})_{\infty},\label{eq:1.25}\end{align}
 they satisfy transformations:\begin{align}
\theta_{1}\left(\frac{v}{\tau}\mid-\frac{1}{\tau}\right) & =-i\sqrt{\frac{\tau}{i}}e^{\pi iv^{2}/\tau}\theta_{1}\left(v\mid\tau\right),\label{eq:1.26}\\
\theta_{2}\left(\frac{v}{\tau}\mid-\frac{1}{\tau}\right) & =\sqrt{\frac{\tau}{i}}e^{\pi iv^{2}/\tau}\theta_{4}\left(v\mid\tau\right),\label{eq:1.27}\\
\theta_{3}\left(\frac{v}{\tau}\mid-\frac{1}{\tau}\right) & =\sqrt{\frac{\tau}{i}}e^{\pi iv^{2}/\tau}\theta_{3}\left(v\mid\tau\right),\label{eq:1.28}\\
\theta_{4}\left(\frac{v}{\tau}\mid-\frac{1}{\tau}\right) & =\sqrt{\frac{\tau}{i}}e^{\pi iv^{2}/\tau}\theta_{2}\left(v\mid\tau\right).\label{eq:1.29}\end{align}
 The Dedekind $\eta(\tau)$ is defined as \cite{Rademarcher}\begin{equation}
\eta(\tau):=e^{\pi i\tau/12}\prod_{k=1}^{\infty}(1-e^{2\pi ik\tau}),\label{eq:1.30}\end{equation}
 or\begin{equation}
\eta(\tau)=q^{1/12}(q^{2};q^{2})_{\infty},\quad q=e^{\pi i\tau},\quad\Im(\tau)>0,\label{eq:1.31}\end{equation}
it has the transformation formula\begin{equation}
\eta\left(-\frac{1}{\tau}\right)=\sqrt{\frac{\tau}{i}}\eta(\tau).\label{eq:1.32}\end{equation}

\begin{lem}
\label{lem:2}For \begin{equation}
0<a<1,\quad n\in\mathbb{N},\quad\gamma>0,\label{eq:1.33}\end{equation}
 and\begin{equation}
q=e^{-2\pi\gamma^{-1}n^{-a}},\label{eq:1.34}\end{equation}
we have\begin{equation}
(q;q)_{\infty}=\sqrt{\gamma n^{a}}\exp\left\{ \frac{\pi}{12}\left((\gamma n^{a})^{-1}-\gamma n^{a}\right)\right\} \left\{ 1+\mathcal{O}\left(e^{-2\pi\gamma n^{a}}\right)\right\} ,\label{eq:1.35}\end{equation}
 and\begin{equation}
\frac{1}{(q;q)_{\infty}}=\frac{\exp\left\{ \frac{\pi}{12}\left(\gamma n^{a}-(\gamma n^{a})^{-1}\right)\right\} }{\sqrt{\gamma n^{a}}}\left\{ 1+\mathcal{O}\left(e^{-2\pi\gamma n^{a}}\right)\right\} \label{eq:1.36}\end{equation}
 as $n\to\infty$. 
\end{lem}
\begin{proof}
From formulas \eqref{eq:1.30}, \eqref{eq:1.31} and \eqref{eq:1.32}
we get \begin{align*}
 & (q;q)_{\infty}=\exp\left(\pi\gamma^{-1}n^{-a}/12\right)\eta\left(\gamma^{-1}n^{-a}i\right)\\
 & =\sqrt{\gamma n^{a}}\exp\left(\pi\gamma^{-1}n^{-a}/12\right)\eta(\gamma n^{a}i)\\
 & =\sqrt{\gamma n^{a}}\exp\left(\pi\gamma^{-1}n^{-a}/12-\pi\gamma n^{a}/12\right)\prod_{k=1}^{\infty}(1-e^{-2\pi\gamma kn^{a}})\\
 & =\sqrt{\gamma n^{a}}\exp\left(\pi\gamma^{-1}n^{-a}/12-\pi\gamma n^{a}/12\right)\left\{ 1+\mathcal{O}\left(e^{-2\pi\gamma n^{a}}\right)\right\} ,\end{align*}
 and\[
\frac{1}{(q;q)_{\infty}}=\frac{\exp\left(\pi\gamma n^{a}/12-\pi\gamma^{-1}n^{-a}/12\right)}{\sqrt{\gamma n^{a}}}\left\{ 1+\mathcal{O}\left(e^{-2\pi\gamma n^{a}}\right)\right\} \]
 as $n\to\infty$. 
\end{proof}

\section{Main Results\label{sec:results}}

For $\Re(z)>-\frac{1}{2}$, we write\begin{equation}
\frac{\Gamma_{q}(z+1/2)}{(q;q)_{\infty}(1-q)^{1/2-z}}=\frac{1}{(q^{z+1/2};q)_{\infty}},\label{eq:2.1}\end{equation}
then,\begin{equation}
\Gamma_{q}(z+\frac{1}{2})=\frac{(q;q)_{\infty}}{(1-q)^{z-1/2}}\sum_{k=0}^{\infty}\frac{q^{k(z+1/2)}}{(q;q)_{k}}.\label{eq:2.2}\end{equation}
Formula \eqref{eq:2.1} implies\begin{equation}
\Gamma_{q}(\frac{1}{2}+z)\Gamma_{q}(\frac{1}{2}-z)=\frac{(1-q)(q;q)_{\infty}^{3}}{(q,q^{1/2-z},q^{1/2+z};q)_{\infty}},\label{eq:2.3}\end{equation}
or\begin{equation}
\Gamma_{q}(\frac{1}{2}+z)\Gamma_{q}(\frac{1}{2}-z)=\frac{(1-q)(q;q)_{\infty}^{3}}{\theta_{4}(q^{z};q^{1/2})}.\label{eq:2.4}\end{equation}
Thus,

\begin{equation}
\Gamma_{q}(\frac{1}{2}-z)=\frac{(q;q)_{\infty}^{2}(1-q)^{z+1/2}}{\theta_{4}(q^{z};q^{1/2})}(q^{z+1/2};q)_{\infty}\label{eq:2.5}\end{equation}
or\begin{equation}
\Gamma_{q}(\frac{1}{2}-z)=\frac{(q;q)_{\infty}^{2}(1-q)^{z+1/2}}{\theta_{4}(q^{z};q^{1/2})}\sum_{k=0}^{\infty}\frac{q^{k(k-1)/2}(-q^{z+1/2})^{k}}{(q;q)_{k}}\label{eq:2.6}\end{equation}
 for $\Re(z)>-\frac{1}{2}$.

\subsection{Case $q\to1$ and $z\to\infty$:}

\begin{thm}
\label{thm:2.3}For \begin{equation}
0<a<\frac{1}{2},\quad n\in\mathbb{N},\quad u\in\mathbb{R},\quad q=\exp(-2n^{-a}\pi),\label{eq:2.7}\end{equation}
we have\begin{align}
\frac{1}{\Gamma_{q}\left(\frac{1}{2}-n-n^{a}u\right)} & =\frac{2\exp\left(\pi n^{-a}(n^{a}u+n)^{2}\right)\cos\pi(n^{a}u+n)\left\{ 1+\mathcal{O}\left(e^{-2\pi n^{a}}\right)\right\} }{\sqrt{n^{a}}\exp\left(\pi n^{a}/12+\pi n^{-a}/6\right)\left(1-\exp(-2\pi n^{-a})\right)^{n+n^{a}u+1/2}},\label{eq:2.8}\end{align}
 and\begin{align}
\frac{1}{\Gamma_{q}\left(\frac{1}{2}+n+n^{a}u\right)} & =\frac{\exp(\pi n^{a}/12-\pi n^{-a}/12)\left\{ 1+\mathcal{O}\left(e^{-2\pi n^{a}}\right)\right\} }{\sqrt{n^{a}}(1-e^{-2\pi n^{-a}})^{1/2-n-n^{a}u}}\label{eq:2.9}\end{align}
 as $n\to\infty$, and the big-O term is uniform with respect $u$
for $u\in[0,\infty)$. 
\end{thm}

\subsection{Case $q\to1$ and $z$ fixed:}

\begin{thm}
\label{thm:2.4}Assume that \begin{equation}
q=e^{-2\pi\tau},\quad\tau>0,\quad x\in\mathbb{R}.\label{eq:2.10}\end{equation}
 If \begin{equation}
x>-\frac{1}{2},\quad q>1-\exp(-2^{2x+1}),\label{eq:2.11}\end{equation}
 then\begin{equation}
\Gamma_{q}(\frac{1}{2}+x)=\Gamma(x+1/2)\left\{ 1+\mathcal{O}\left((1-q)\log^{2}(1-q)\right)\right\} ,\label{eq:2.12}\end{equation}
 and\begin{equation}
\frac{1}{\Gamma_{q}(\frac{1}{2}-x)}=\frac{\left\{ 1+\mathcal{O}\left((1-q)\log^{2}(1-q)\right)\right\} }{\Gamma(\frac{1}{2}-x)},\label{eq:2.13}\end{equation}
 where the implicit constants of the big-O terms are independent of
$x$ under the condition \eqref{eq:2.11}
\end{thm}

\section{Proofs\label{sec:proofs}}

\subsection{Proof for Theorem \ref{thm:2.3}}

\begin{proof}
We first observe that\[
\frac{1}{(q;q)_{\infty}}=\frac{\exp\left(\pi n^{a}/12-\pi n^{-a}/12\right)}{\sqrt{n^{a}}}\left\{ 1+\mathcal{O}\left(e^{-2\pi n^{a}}\right)\right\} ,\]
and \begin{align*}
 & \frac{1}{\Gamma_{q}\left(1/2-n-n^{a}u\right)}=\frac{(q^{1/2-n}e^{2\pi u};q)_{\infty}}{(q;q)_{\infty}(1-q)^{n+n^{a}u+1/2}}\\
 & =\frac{(q,q^{1/2}e^{-2\pi u},q^{1/2}e^{2\pi u};q)_{\infty}q^{-n^{2}/2}e^{2\pi nu}}{(-1)^{n}(1-q)^{n+n^{a}u+1/2}(q,q,q^{n+1/2}e^{-2\pi u};q)_{\infty}}\end{align*}
as $n\to\infty$. Then we have\begin{align*}
\frac{1}{(q,q,q^{n+1/2}e^{-2\pi u};q)_{\infty}} & =n^{-a}\exp\left(\pi n^{a}/6-\pi n^{-a}/6\right)\left\{ 1+\mathcal{O}\left(e^{-2\pi n^{a}}\right)\right\} ,\end{align*}
and \begin{align*}
 & (q,q^{1/2}e^{-2\pi u},q^{1/2}e^{2\pi u};q)_{\infty}=\theta_{4}(ui\mid n^{-a}i)=n^{a/2}e^{\pi n^{a}u^{2}}\theta_{2}(n^{a}u\mid n^{a}i)\\
 & =2n^{a/2}\exp\pi n^{a}(u^{2}-1/4)\cos(n^{a}u\pi)\left\{ 1+\mathcal{O}\left(e^{-2\pi n^{a}}\right)\right\} \end{align*}
 as $n\to\infty$. Thus,\begin{align*}
\frac{1}{\Gamma_{q}\left(1/2-n-n^{a}u\right)} & =\frac{2\exp\pi\left(n^{-a}(n^{a}u+n)^{2}\right)\cos(\pi n^{a}u+n\pi)\left\{ 1+\mathcal{O}\left(e^{-2\pi n^{a}}\right)\right\} }{n^{a/2}\left(1-e^{-2\pi n^{-a}}\right)^{n+n^{a}u+1/2}\exp\left(\pi n^{a}/12+\pi n^{-a}/6\right)}\end{align*}
 as $n\to\infty$, and it is clear that the big-O term is uniform
with respect to $u\ge0$.

Similarly, formula \eqref{eq:2.9} follows from Lemma \ref{lem:1}
and \ref{lem:2}. 
\end{proof}

\subsection{Proof for Theorem \ref{thm:2.4}}

\begin{proof}
From \eqref{eq:1.27}, \eqref{eq:1.32} and\[
\Gamma_{q}(\frac{1}{2}+x)\Gamma_{q}(\frac{1}{2}-x)=\frac{\eta(\tau i)^{3}e^{\pi\tau/4}(1-e^{-2\pi\tau})}{\theta_{4}(x\tau i\vert\tau i)},\]
we get\[
\Gamma_{q}(\frac{1}{2}+x)=\frac{(1-e^{-2\pi\tau})\exp\left(-\pi\tau(x^{2}-1/4)\right)}{\tau\Gamma_{q}(\frac{1}{2}-x)}\frac{\eta^{3}(i/\tau)}{\theta_{2}(x\vert i/\tau)},\]
 and \eqref{eq:1.8} and \eqref{eq:1.15} imply that\begin{align*}
 & \Gamma_{q}(\frac{1}{2}+x)=\frac{(1-e^{-2\pi\tau})\exp\left(-\pi\tau(x^{2}-1/4)\right)}{2\tau\cos(\pi x)\Gamma_{q}(\frac{1}{2}-x)}\left\{ 1+\mathcal{O}\left(\exp(-\frac{2\pi}{\tau}\right)\right\} \\
 & =\frac{\pi\exp\left(-\pi\tau x^{2}\right)\left\{ 1+\mathcal{O}\left(\tau\right)\right\} }{\cos(\pi x)\Gamma_{q}(\frac{1}{2}-x)}\\
 & =\frac{\pi\exp\left(-\pi\tau\log^{2}(1-q)\frac{x^{2}}{\log^{2}(1-q)}\right)\left\{ 1+\mathcal{O}\left(\tau\right)\right\} }{\cos(\pi x)\Gamma_{q}(\frac{1}{2}-x)}\end{align*}
 as $\tau\to0^{+}$ and the big-O term is independent of $x$.

The condition \eqref{eq:2.11} implies that\[
0<\frac{x^{2}}{\log^{2}(1-q)}<1,\]
then,\[
\Gamma_{q}(\frac{1}{2}+x)=\frac{\pi\left\{ 1+\mathcal{O}\left((1-q)\log^{2}(1-q)\right)\right\} }{\cos(\pi x)\Gamma_{q}(\frac{1}{2}-x)}\]
 as $q\to1$ and the implicit constant above is independent of $x$.

It is well-known that an $q$-analogue of \eqref{eq:1.6} is \cite{Andrews4,Gasper,Ismail2,Koekoek}
\[
\int_{0}^{\infty}\frac{t^{x-1/2}}{(-t;q)_{\infty}}dt=\frac{\pi}{\cos\pi x}\frac{(q^{1/2-x};q)_{\infty}}{(q;q)_{\infty}},\quad x>-\frac{1}{2},\]
 or\[
\int_{0}^{\infty}\frac{t^{x-1/2}dt}{(-(1-q)t;q)_{\infty}}=\frac{\pi}{\cos(\pi x)\Gamma_{q}(\frac{1}{2}-x)},\quad x>-\frac{1}{2}.\]
 Consequently,\[
\Gamma_{q}(\frac{1}{2}+x)=\left\{ 1+\mathcal{O}\left((1-q)\log^{2}(1-q)\right)\right\} \int_{0}^{\infty}\frac{t^{x-1/2}dt}{(-(1-q)t;q)_{\infty}}\]
 as $\tau\to0^{+}$ and the implicit constant of the big-O term here
is independent of $x$ with $x>-\frac{1}{2}$. 

Write \[
\int_{0}^{\infty}\frac{t^{x-1/2}dt}{(-(1-q)t;q)_{\infty}}:=I_{1}+I_{2},\]
 where\[
I_{1}:=\int_{0}^{\log(1-q)^{-2}}\frac{t^{x-1/2}dt}{(-(1-q)t;q)_{\infty}},\]
 and\[
I_{2}:=\int_{\log(1-q)^{-2}}^{\infty}\frac{t^{x-1/2}dt}{(-(1-q)t;q)_{\infty}}.\]
 In $I_{1}$ we have\begin{eqnarray*}
\log(-(1-q)t;q)_{\infty} & = & \sum_{k=0}^{\infty}\log(1+(1-q)tq^{k})\\
 & = & \sum_{k=0}^{\infty}\sum_{n=0}^{\infty}\frac{(-1)^{n}}{n+1}\left((1-q)tq^{k}\right)^{n+1}\\
 & = & \sum_{n=0}^{\infty}\frac{(-1)^{n}t^{n+1}}{n+1}\frac{(1-q)^{n+1}}{1-q^{n+1}}\\
 & = & t+r(n),\end{eqnarray*}
 where\[
r(n):=\sum_{n=1}^{\infty}\frac{(-1)^{n}t^{n+1}}{n+1}\frac{(1-q)^{n+1}}{1-q^{n+1}}.\]
The condition \eqref{eq:2.11} implies that\[
0<(1-q)\log(1-q)^{-1}<e^{-1},\]
then\begin{eqnarray*}
|r(n)| & \le & c(1-q)\log^{2}(1-q),\end{eqnarray*}
with\[
c=4\sum_{n=0}^{\infty}\frac{\left[2(1-q)\log(1-q)^{-1}\right]^{n}}{n+2}<4\sum_{n=0}^{\infty}\frac{(2e^{-1})^{n}}{n+2}.\]
Therefore,\[
I_{1}=\int_{0}^{2\log(1-q)^{-1}}e^{-t}t^{x-1/2}dt\left\{ 1+\mathcal{O}\left((1-q)\log^{2}(1-q)\right)\right\} \]
as $\tau\to0^{+}$, and the implicit constant of the big-O term is
independent of $x$. 

Since \begin{align*}
\int_{2\log(1-q)^{-1}}^{\infty}e^{-t}t^{x-1/2}dt & \le e^{-\log(1-q)^{-1}}\int_{2\log(1-q)^{-1}}^{\infty}e^{-t/2}t^{x-1/2}dt<(1-q)\Gamma(x+1/2)2^{x+1/2},\end{align*}
 clearly, \[
\int_{2\log(1-q)^{-1}}^{\infty}e^{-t}t^{x-1/2}dt=\Gamma(x+1/2)\mathcal{O}\left((1-q)\log^{2}(1-q)\right).\]
 Thus,\[
I_{1}=\Gamma(x+1/2)\left\{ 1+\mathcal{O}\left((1-q)\log^{2}(1-q)\right)\right\} ,\]
 as $\tau\to0^{+}$, and the implicit constant of the big-O term is
independent of $x$ under the condition \eqref{eq:2.11}. 

Recall that \[
(-(1-q)t;q)_{\infty}=\sum_{n=0}^{\infty}q^{n(n-1)/2}t^{n}\frac{(1-q)^{n}}{(q;q)_{n}}>q^{n(n-1)/2}t^{n}\frac{(1-q)^{n}}{(q;q)_{n}},\]
 for any $n=\left\lfloor -\log(1-q)\right\rfloor \ge2^{2x+1}>2x+1$.
Then,\begin{align*}
I_{2} & \le\frac{q^{n(1-n)/2}(q;q)_{n}(-\log(1-q))^{2x-2n+1}}{(1-q)^{n}(n-x-1/2)}\\
 & <\frac{\Gamma(x+1/2)n!q^{n(1-n)/2}(-\log(1-q))^{2x-2n+1}}{\Gamma(x+3/2)}\\
 & <\Gamma(x+1/2)n!q^{n(1-n)/2}(-\log(1-q))^{2x-2n+1}.\end{align*}
 It is clear that\[
q^{n(1-n)/2}=\mathcal{O}(1)\]
 as $\tau\to0^{+}$. From the Stirling formula\[
n!=\sqrt{2\pi n}\left(\frac{n}{e}\right)^{n}\left\{ 1+\mathcal{O}\left(\frac{1}{n}\right)\right\} \]
as $n\to\infty$, we have\[
I_{2}=\Gamma(x+\frac{1}{2})\mathcal{O}\left((1-q)\log^{1/2}(1-q)^{-1}\right)\]
as $\tau\to0^{+}$ and the implicit constant of the big-O term is
independent of $x$ under the condition \eqref{eq:2.11}.

Therefore, under the condition \eqref{eq:2.11}\[
\int_{0}^{\infty}\frac{t^{x-1/2}dt}{(-(1-q)t;q)_{\infty}}=\Gamma(x+1/2)\left\{ 1+\mathcal{O}\left((1-q)\log^{2}(1-q)\right)\right\} \]
as $\tau\to0^{+}$ and the implicit constant of the big-O term is
independent of $x>-\frac{1}{2}$. 

Hence we have proved that under the condition \eqref{eq:2.11} we
have\[
\Gamma_{q}(\frac{1}{2}+x)=\Gamma(x+\frac{1}{2})\left\{ 1+\mathcal{O}\left((1-q)\log^{2}(1-q)\right)\right\} ,\quad x>-1/2.\]
 Then,\[
\Gamma_{q}(\frac{1}{2}-x)=\frac{\pi\sec\pi x}{\Gamma(x+1/2)}\left\{ 1+\mathcal{O}\left((1-q)\log^{2}(1-q)\right)\right\} ,\]
which is\[
\frac{1}{\Gamma_{q}(\frac{1}{2}-x)}=\frac{\left\{ 1+\mathcal{O}\left((1-q)\log^{2}(1-q)\right)\right\} }{\Gamma(\frac{1}{2}-x)},\]
where the implicit constant of the big-O term is independent of $x>-\frac{1}{2}$. 
\end{proof}

\end{document}